\documentclass[11pt]{article}
\usepackage{amsmath}
\usepackage{cases}
\usepackage{amsfonts,latexsym,amsmath,amscd,geometry}
\usepackage{amssymb}
\usepackage{latexsym}
\usepackage{hyperref, xcolor}

\newcommand \nc{\newcommand}
\newtheorem{theorem}{Theorem}[section]
\newtheorem{lemma}[theorem]{Lemma}

\newtheorem{corollary}[theorem]{Corollary}

\newtheorem{remark}[theorem]{Remark}

\nc{\ba}{\begin{array}}\nc{\ea}{\end{array}}
\nc{\be}{\begin{eqnarray}}\nc{\ee}{\end{eqnarray}}
\nc{\beq}{\begin{equation}}\nc{\eeq}{\end{equation}}
\nc{\bex}{\begin{eqnarray*}}\nc{\eex}{\end{eqnarray*}}
\nc{\btm}{\begin{theorem}} \nc{\etm}{\end{theorem}}
\nc{\blm}{\begin{lemma}} \nc{\elm}{\end{lemma}}
\nc{\R}{\mathbb{R}} \nc{\va}{\varepsilon} \nc{\ls}{\limits}

\def\pf{\noindent{\bf Proof.\quad}}\def\endpf{\hfill$\Box$}

\begin{document}
\title{Regularity and uniqueness for a class of solutions to the hydrodynamic flow of nematic liquid crystals}
\author{Tao Huang \\
Department of Mathematics, The Pennsylvania State University\\
University Park, PA 16802, USA\\
txh35@psu.edu}
\date{}
\maketitle

\begin{abstract} In this paper, we establish an $\epsilon$-regularity criterion for any weak solution
$(u,d)$ to the nematic liquid crystal flow (\ref{lce}) such that
$(u,\nabla d)\in L^p_tL^q_x$ for some $p\ge 2$ and $q\ge n$ satisfying the  condition
(\ref{serrin-condition}). As consequences, we prove the interior smoothness of any such a solution
when $p>2$ and $q>n$. We also show that uniqueness holds for the class of weak solutions $(u,d)$
the Cauchy problem of the nematic liquid crystal flow (\ref{lce}) that satisfy $(u,\nabla d)\in L^p_tL^q_x$
for some $p>2$ and $q>n$ satisfying (\ref{serrin-condition}).
\end{abstract}

\section {Introduction}
\setcounter{equation}{0}
\setcounter{theorem}{0}
For any $n\geq 3$,
the hydrodynamic flow of nematic liquid crystals in $\mathbb R^n\times[0,T]$, for some $0<T<+\infty$, is given by
\begin{equation}\label{lce}
\begin{cases}
u_t+u\cdot\nabla u-\Delta u+\nabla P=-\nabla\cdot(\nabla d\otimes\nabla d-\frac12{|\nabla d|^2}\mathbb{I}_n)
 & \ {\rm{in}}\ \mathbb R^n\times (0,T)\\
\nabla\cdot u=0 & \ {\rm{in}}\ \mathbb R^n\times (0,T)\\
d_t+u\cdot\nabla d=\Delta d+|\nabla d|^2d & \ {\rm{in}}\ \mathbb R^n\times (0,T)\\
(u,d)=(u_0,d_0)& \ {\rm{on}}\ \mathbb R^n\times\{0\}
\end{cases}
\end{equation}
where $u:\mathbb R^n\times[0,T]\rightarrow\R^n$ is  the velocity field of underlying incompressible fluid,
$d:\mathbb R^n\times[0,T]\rightarrow S^2$ is the director field of nematic liquid crystal molecules,
$P:\mathbb R^n\times[0,T]\rightarrow\R$ is the pressure function,
$\nabla\cdot$ denotes the divergence operator on $\mathbb R^n$, $\nabla
d\otimes\nabla d =\left(\frac{\partial d}{\partial
x_i}\cdot\frac{\partial d}{\partial x_j}\right)_{1\leq i,j\leq n}$
is the stress tensor induced by the director field $d$,
$\mathbb{I}_n$ is the identity matrix of order $n$, $u_0:\mathbb R^n\to\mathbb R^n$ is the initial
velocity field with $\nabla\cdot u_0=0$, and $d_0:\mathbb R^n\to S^2$ is the initial
director field.

The system (\ref{lce}) is a simplified version of the
Ericksen-Leslie system modeling the hydrodynamics of liquid crystal
materials, proposed by Ericksen \cite{ericksen} and Leslie
\cite{leslie} in 1960's. It is a macroscopic continuum description
of the time evolution of the material under the influence of both
the flow field and the macroscopic description of the microscopic
orientation configurations of rod-like liquid crystals.
 The interested readers can refer to
\cite{ericksen}, \cite{leslie}, \cite{lin}, and \cite{lin-liu} for
more detail. Mathematically, the system (\ref{lce}) is strongly
coupling the Naiver-Stokes equations and the
(transported) heat flow of harmonic maps into $S^2$.

For $n=2$, Lin-Lin-Wang \cite{lin-lin-wang} have proved the existence of global Leray-Hopf type weak solutions to (\ref{lce}) with initial and boundary conditions, which is smooth away from finitely many possible singular times
(see Hong \cite{hong} and Xu-Zhang \cite{xu-zhang} for related works). Lin-Wang \cite{lin-wang} proved the uniqueness for such weak solutions.  It remains a very challenge open problem to prove the global existence of Leray-Hopf type weak solutions
and partial regularity of suitable weak solutions to (\ref{lce}) in higher dimensions. A BKM type blow-up criterion was obtained
for the local strong solution to (\ref{lce}) for $n=3$ by \cite{huang-wang2}, i.e., if $0<T_*<+\infty$ is the maximum time
interval of the strong solution to (\ref{lce}), then
$$\int_0^{T_*}\left(\|\nabla\times u\|_{L^{\infty}}+\|\nabla d\|^2_{L^{\infty}}\right)\,dt=+\infty.$$
Recently, the local well-posedness of (\ref{lce}) was obtained for initial data $(u_0,d_0)$ with $(u_0,\nabla d_0)\in L^3_{\rm{uloc}}(\mathbb R^3)$, the space of uniformly locally $L^3$-integrable functions,  of small norm for $n=3$
by \cite{hineman-wang}.  While the global well-posedness of (\ref{lce}) was obtained by \cite{wang1}
for $(u_0,d_0)\in$ BMO$\times$ BMO$^{-1}$ of small norm for $n\ge 3$.

The existence of global Leray-Hopf weak
solutions to the Naiver-Stokes equations has long been established by Leray \cite{leray} and
Hopf \cite{hopf}. However the uniqueness (regularity) of Leray-Hopf solutions
in dimension three remains largely open. In
\cite{serrin}, Serrin proved the so called `weak-strong' uniqueness,
i.e., the uniqueness holds for Leray-Hopf solutions $u,v$ with the same initial data, if
$u\in L^{p}_tL^q_x(\mathbb R^n\times[0,T])$, where $p\geq 2$ and $q\ge n$ satisfy
\begin{equation}\label{serrin-condition}
\frac{2}{p}+\frac{n}{q}=1.
\end{equation}
The smoothness of such solutions was established by Ladyzhenskaya in \cite{ladyzhen} for $p>2$ and $q>n$.
In the fundamental work \cite{ESS}, Escauriaza-Seregin-$\check{\mbox{S}}$ver$\acute{\mbox{a}}$k have proved the smoothness of Serrin's solutions for the endpoint case $(p,q)=(+\infty,n)$ when $n=3$ (see also \cite{dong-du} for $n\ge 4$).
Wang \cite{wang} proved smoothness of weak solutions $u$ to the heat flow of harmonic maps
such that $\nabla u\in L^{p}_tL^q_x(\mathbb R^n\times[0,T])$ with $\frac{2}{p}+\frac{n}q=1$
for $n\geq 4$ (or $q\geq 4$ for $2\leq n<4$, see \cite{huang-wang1} for the case $2<q<4$ when $2\le n<4$).
In \cite{huang-wang1},  the uniqueness of Serrin's solutions to the heat flow of harmonic maps is also established
when $p> 2$ and $q>n$. These results motivate us
to investigate the regularity and uniqueness of Serrin's ($p,q$)-solutions to the system (\ref{lce}) of
nematic liquid crystal flows.

Before stating our main theorems, we need to introduce some notations.

\medskip
\noindent{\bf Notations}: For $1\le p, q\le +\infty$, $0<T\le \infty$, define the Sobolev space
$$H^{1}(\mathbb R^n\times[0,T], \mathbb R^n)
=\Big\{f\in L^2([0,T], H^{1}(\mathbb R^n, \mathbb R^n)): \ \partial_t
f\in L^2([0,T], L^2(\mathbb R^n,\mathbb R^n))\Big\},$$
$$\mathbb E^p(\mathbb R^n\times[0,T], \mathbb R^n)=\{f\in L^p(\mathbb R^n\times[0,T], \mathbb R^n)
\ | \ \nabla\cdot f=0\},$$
and the Morrey space $M^{p,\lambda}(U)$ for $0\le\lambda\le n+2$
and $U=U_1\times [t_1,t_2]\subset\mathbb R^{n}\times\mathbb R$:
$${M}^{p,\lambda}(U)=\left\{f\in L^{p}_{\mbox{loc}}(U):
\Big\|f\Big\|_{{M}^{p,\lambda}(U)}<+\infty\right\},$$
where
$$\Big\|f\Big\|_{{M}^{p,\lambda}(U)}=\Big(\sup\limits_{(x,t)\in U}\sup\limits_{0<r<\min\delta((x,t),\partial_p U)}\
r^{\lambda-n-2}\int_{P_r(x,t)} |f|^p\Big)^{\frac{1}{p}},$$
$$B_r(x)=\{y\in\mathbb R^n: \,|y-x|\leq r\},
\ P_r(x,t)=B_r(x)\times[t-r^2,t], \partial_p U=(\partial U_1\times [t_1,t_2])\cup (U_1\times\{t_1\}),$$
and
$$\delta((x,t),\partial_p U)=\inf_{(y,s)\in\partial_p U}\delta\left((x,t), (y,s)\right),
\ {\rm{and}}\  \  \delta((x,t), (y,s))=\min\left\{|x-y|, \sqrt{|t-s|}\right\}.$$
Denote $B_r$ (or $P_r$) for
$B_r(0)$ (or $P_r(0)$ respectively). We also
recall the weak Morrey space,  $M_*^{p,\lambda}(U)$,
that is the set of functions $f$ on $U$ such that
$$\|f\|^p_{M_*^{p,\lambda}(U)}=\sup\limits_{r>0,(x,t)\in U }\,\Big\{r^{\lambda-(n+2)}
\|f\|_{L^{p,*}(P_r(x,t)\cap U)}^p \Big\}<+\infty,$$
where $L^{p,*}(P_r(x,t)\cap U)$ is the weak $L^p$-space, that is the collection of functions $v$ on $P_r(x,t)\cap U$
such that
$$
\|v\|^p_{L^{p,*}(P_r(x,t)\cap U)}=\sup\limits_{a>0}\,  \Big\{\ a^p\left|\{z\in P_r(x,t)\cap U\ : \, |v(z)|>a\}\right|\Big\}<+\infty.
$$

Recall that $(u,d)\in H^{1}(\mathbb R^n\times[0,T],\mathbb R^n\times S^2)$ is a weak solution to
(\ref{lce}) if $(u,d)$ satisfies
(\ref{lce})$_1$-(\ref{lce})$_3$ in sense of distribution and (\ref{lce})$_4$ in sense of trace.
A weak solution $(u,d)\in H^{1}(\mathbb R^n\times[0,T],\mathbb R^n\times S^2)$ of (\ref{lce}) if called
a Serrin's ($p,q$)-solution, if $(u, \nabla d)\in L^p_tL^q_x(\mathbb R^n\times [0,T])$ for some
($p,q$) satisfying (\ref{serrin-condition}).
Our first result concerns an $\epsilon_0$-regularity criterion for Serrin's ($p,q$)-solutions to (\ref{lce}).

\btm\label{loc-reg-th-lcf}{\it There exists $\epsilon_0>0$ such that
if a weak solution $(u,d)\in H^1(P_r(x_0,t_0), \mathbb R^n\times S^2)$ to (\ref{lce}) satisfies
 \beq\label{u-lcf2.1}
\|u\|_{L^{p}_t L^{q}_x(P_r(x_0,t_0))}+\|\nabla
d\|_{L^{p}_t L^{q}_x(P_r(x_0,t_0))}\leq\epsilon_0,
 \eeq where $p\geq 2$ and
$q\geq n$ satisfy (\ref{serrin-condition}),
then $(u,d)\in
C^{\infty}(P_{\frac{r}{16}}(x_0,t_0))$,  and
 \beq\label{u-lcf2.2}
r\|u\|_{L^{\infty}(P_{\frac{r}{16}}(x_0,t_0))}+r\|\nabla
d\|_{L^{\infty}(P_{\frac{r}{16}}(x_0,t_0))}
\leq C\left(\|u\|_{L^{p}_t L^{q}_x(P_r(x_0,t_0))}+\|\nabla
d\|_{L^{p}_t L^{q}_x(P_r(x_0,t_0))}\right).
\eeq
}
 \etm

 A direct corollary of Theorem \ref{loc-reg-th-lcf} is the following regularity theorem
for Serrin's ($p,q$)-solutions to (\ref{lce}).
\begin{corollary}\label{reg-th-lcf}{\it For some $0<T< +\infty$,
suppose that $(u,d)\in H^1(\mathbb R^n\times[0,T],\mathbb R^n\times S^2)$ is a weak solution to (\ref{lce}) with $(u,\nabla d)\in L^p_tL^q_x(\mathbb R^n\times[0,T])$, for some $p>2$ and
$q>n$ satisfying (\ref{serrin-condition}). Then $(u,d)\in
C^{\infty}(\mathbb R^n\times(0,T],\mathbb R^n\times S^2)$.
}
\end{corollary}

\begin{remark}{\it (i) For the heat flow of harmonic maps and the Navier-Stokes equations, Corollary \ref{reg-th-lcf} is valid for the end point case $(p,q)=(+\infty,n)$. It is an interesting open question to investigate the regularity of Serrin's solutions
to (\ref{lce})
in this end point case.  \\
(ii) If $(u_0,\nabla d_0)\in L^{\gamma}(\mathbb R^n)$ for some $\gamma>n$, then the local existence of Serrin's solutions in $L^p_tL^q_x$ for some $p>2$ and $q>n$ can be obtained by the fixed point argument (see e.g., \cite{FJR} $\S 4$).\\


}
\end{remark}

As a corollary of Theorem \ref{loc-reg-th-lcf} and Corollary \ref{reg-th-lcf}, we can prove the uniqueness of
Serrin's ($p,q$)-solutions to (\ref{lce}).

\btm\label{unique3}{\it For $n\geq 2$, $0<T< +\infty$, and
$i=1,2$,  if $(u_i, d_i):\mathbb R^n\times[0,T]\to\mathbb R^n\times S^2$ are two
weak solutions to (\ref{lce}) with the same initial
data $(u_0,d_0):\mathbb R^n\to \mathbb R^n\times S^2$. Suppose, in
additions, there exists $p>2$ and $q>n$ satisfying
(\ref{serrin-condition}) such that  $(u_1,\nabla d_1), (u_2,\nabla
d_2)\in L^{p}_tL^{q}_x(\mathbb R^n\times[0,T])$. Then $(u_1,d_1)\equiv (u_2,d_2)$ on
$\mathbb R^n\times[0,T]$.
 }\etm

\begin{remark} {\rm For $n=2$,  Lin-Wang \cite{lin-wang} have proved the uniqueness of (\ref{lce})
for $p=q=4$.  More precisely, if, for $i=1,2$,
$$
\begin{cases}
u_i\in L^{\infty}([0,T], L^{2}(\R^2,\R^2))
\cap L^{2}([0,T], H^{1}(\R^2,\R^2));\\
\nabla d_i\in L^{\infty}([0,T], L^2(\R^2))\cap L^{2}([0,T], {H}^{1}(\R^2))
\end{cases}
$$
are weak solutions to (\ref{lce}) under the same initial condition,
then $(u_1,d_1)\equiv (u_2,d_2)$ on $\R^2\times[0,T]$. For $n\geq 3$, Lin-Wang \cite{lin-wang}
proved the uniqueness,
provided that $u_i\in C([0,T], L^{n}(\R^n))$ and $\nabla d_i\in C([0,T], {L}^{n}(\R^n))$ for $i=1,2$.}

\end{remark}

\section{Proof of Theorem \ref{loc-reg-th-lcf} and Corollary \ref{reg-th-lcf}}
\setcounter{equation}{0}
\setcounter{theorem}{0}
In this section, we will prove Theorem
\ref{loc-reg-th-lcf} and Corollary \ref{reg-th-lcf} for nematic liquid crystal flows (\ref{lce}). The crucial step
is to establish an $\epsilon_0$-regularity criterion.

\blm\label{unqlc-lemma3.1}{\it There exists $\epsilon_0>0$ such that
if $(u,\nabla d)\in L^p_tL^q_x(P_1(0,1))$, for some $p\geq 2$ and
$q\geq n$ satisfying (\ref{serrin-condition}), is a weak
solution to (\ref{lce}) that satisfies
 \beq\label{unqlc3.3}
\|u\|_{L^{p}_t L^{q}_x(P_1(0,1))}+\|\nabla
d\|_{L^{p}_t L^{q}_x(P_1(0,1)))}\leq\epsilon_0,
 \eeq
then $(u,d)\in
C^{\infty}(P_{\frac{1}{16}}(0,1))$,  and
 \beq\label{unqlc3.4}
\|u\|_{L^{\infty}(P_{\frac{1}{16}}(0,1))}+\|\nabla
d\|_{L^{\infty}(P_{\frac{1}{16}}(0,1))}\leq C\epsilon_0.
\eeq
}\elm

Before proving this lemma, we need the following inequality, due to Serrin \cite{serrin}.

\blm\label{unq-serrin-lemma}{\it For any open set $U\subset\mathbb
R^n$ and any open interval $I\subset\mathbb R$, let $f$, $g$, $h\in
L^2_tH^1_x(U\times I)$ and $f\in L^p_tL^q_x(U\times I)$ with $3\le
n\le q\le +\infty$ and $2\le p\le +\infty$ satisfying
(\ref{serrin-condition}). Then
\beq\label{unq2.5}
\int_{U\times I}|f||g||\nabla h| \leq C\|\nabla h\|_{L^2(U\times I)}
\|g\|_{L^2_tH^1_x(U\times I)}^{\frac{n}{q}}
\left\{\int_I\|f\|_{L^q(\mathbb R^n)}^{p}\|g\|_{L^2(\mathbb R^n)}^2\,dt\right\}^{\frac{1}{p}},
\eeq where $C>0$ depends only on $n$. }\elm

\noindent{\bf Proof of Lemma \ref{unqlc-lemma3.1}}.
For any $(x,t)\in
P_{\frac{1}{2}}(0,1)$ and $0<r<\frac{1}{2}$, we have, by
(\ref{unqlc3.3}),
 \beq\label{unqlc3.5}
\|u\|_{L^{p}_tL^q_x(P_r(x,t))}+\|\nabla
d\|_{L^{p}_tL^q_x(P_r(x,t))}\leq\epsilon_0.
\eeq
We will divide the proof into two claims.

\noindent{\bf Claim 1}. $\nabla d\in L^{\gamma}(P_{\frac{1}{2}}(0,1))$ for any $1<\gamma<\infty$, and
\begin{equation}
\|\nabla d\|_{L^{\gamma}(P_{\frac{1}{4}}(0,1))}\leq C(\gamma)\|\nabla d\|_{L^p_tL^q_x(P_1(0,1))}.
\end{equation}
To show it,  let $d_1: P_r(x,t)\rightarrow \R^3$ solve
\beq\label{unqlc3.6}
\left\{
\begin{split}
\partial_t d_1-\Delta d_1&=0, \quad\mbox{in } P_r(x,t)\\
d_1&=d, \quad\mbox{on }\partial_pP_r(x,t).
\end{split}\right.
\eeq Set $d_2=d-d_1$. Multiplying (\ref{lce})$_3$ and
(\ref{unqlc3.6}) by $d_2$, subtracting the resulting equations and
integrating over $P_r(x,t)$, we obtain
 \beq\label{unqlc3.7}
\begin{split}
&\sup\limits_{t-r^2\leq \tau\leq t}\int_{B_r(x)}|d_2|^2(\cdot,\tau)+2\int_{P_r(x,t)}|\nabla d_2|^2\\
\leq& C\int_{P_r(x,t)}(|u||d_2||\nabla d|+|\nabla d||d_2||\nabla d|)=J_1+J_2.
\end{split}
\eeq By (\ref{unq2.5}), the Poincar$\acute{\mbox{e}}$ inequality and
the Young inequality, we have
 \begin{eqnarray*}
 |J_1|
 &\lesssim&\begin{cases}
 \|\nabla d\|_{L^2(P_r(x,t))} \|\nabla d_2\|_{L^2(P_r(x,t))}^{\frac{n}{q}}
 \left\{\int_{t-r^2}^t\|u\|_{L^{q}(B_{r}(x))}^{p}\|d_2\|_{L^2(B_r(x))}^2\,d\tau\right\}^{\frac{1}{p}},\ & p<+\infty\\
 \|\nabla d\|_{L^2(P_r(x,t))}\|\nabla d_2\|_{L^2(P_r(x,t))}
 \|u\|_{L^{\infty}_tL^n_x(P_r(x,t))},\ & p=+\infty,
\end{cases}\\
 &\leq& \begin{cases}\frac{1}{2}\|\nabla
 d_2\|_{L^2(P_r(x,t))}^{2}+C\epsilon_0\|\nabla
 d\|_{L^2(P_r(x,t))}^2
 +C\epsilon_0^{\frac{p}2}\|d_2\|_{L^\infty_tL^2_x(P_r(x,t))}^2,\ & p<+\infty\\
\frac{1}{2}\|\nabla
 d_2\|_{L^2(P_r(x,t))}^{2}+C\epsilon_0\|\nabla
 d\|_{L^2(P_r(x,t))}^2, \ & p=+\infty.
\end{cases}
\end{eqnarray*}
Similarly, for $J_2$, we have
$$
|J_2| \leq \begin{cases}\frac{1}{2}\|\nabla
 d_2\|_{L^2(P_r(x,t))}^{2}+C\epsilon_0\|\nabla
 d\|_{L^2(P_r(x,t))}^2
 +C\epsilon_0^{\frac{p}2}\|d_2\|_{L^\infty_tL^2_x(P_r(x,t))}^2,\ & p<+\infty\\
\frac{1}{2}\|\nabla
 d_2\|_{L^2(P_r(x,t))}^{2}+C\epsilon_0\|\nabla
 d\|_{L^2(P_r(x,t))}^2, \ & p=+\infty.
\end{cases}
$$
Putting these estimates into (\ref{unqlc3.7}), applying (\ref{unqlc3.5}), and choosing sufficiently small $\epsilon_0$,
 we have
\beq\label{unqlc3.12}
\int_{P_r(x,t)}|\nabla d_2|^2
\leq C \epsilon_0\|\nabla d\|_{L^2(P_r(x,t))}^2.
\eeq
This, combined with the standard estimate on $d_1$, implies that for any $\theta\in (0,1)$,
\beq\label{unqlc3.14}
(\theta r)^{-n}\int_{P_{\theta r}(x,t)}|\nabla d|^2\le
 C\Big(\theta^2+\theta^{-n}\epsilon_0\Big) r^{-n}\int_{P_r(x,t)}|\nabla d|^2.
\eeq
By iterations, we obtain for any $(x,t)\in P_{\frac{1}{2}}(0,1)$,
$0<r\leq\frac{1}{2}$ and $0<\alpha<1$,
\beq\label{unqlc3.15}
\begin{split}
r^{-n}\int_{P_{r}(x,t)}|\nabla d|^2
\leq Cr^{2\alpha}\int_{P_1(0,1)}|\nabla d|^2.
\end{split}
\eeq
Hence $\nabla d\in\mathcal{M}^{2,2-2\alpha}(P_\frac12(0,1))$ and
\begin{equation}\label{Dd_morrey_estimate}
\|\nabla d\|_{\mathcal M^{2,2-2\alpha}(P_\frac12(0,1))}\le
C\|\nabla d\|_{L^p_tL^q_x(P_1(0,1))}.
\end{equation}
Now Claim 1 follows
by the same estimate of Riesz potentials between parabolic Morrey spaces
as in \cite{huang-wang} (Theorem 1.5) and \cite{lin-wang} (Lemma 2.1).

\medskip
\noindent{\bf Claim 2}. $u\in L^{\gamma}(P_{\frac{1}{4}}(0,1))$ for
any $1<\gamma<\infty$, and
\begin{equation} \label{integral_estimate_u}
\|u\|_{L^\gamma(P_\frac14(0,1))}\le C(\gamma)\|u\|_{L^p_tL^q_x(P_1(0,1))}.
\end{equation}

Let $\mathbb{E}^\gamma$ be the closure in
$L^{\gamma}(\R^n,\R^n)$ of all divergence-free vector fields with
compact supports. Let $\mathbb{P}:L^{2}(\R^n,\R^n)\rightarrow
\mathbb{E}^2$ be the Leray projection operator. It is well-known
that $\mathbb{P}$ can be extended to a bounded linear operator from
$L^\gamma(\R^n,\R^n)$ to $\mathbb{E}^{\gamma}$ for all
$1<\gamma<+\infty$. Let $\mathbb{A}=\mathbb{P}\Delta$ denote the Stokes
operator.

For any $(x,t)\in P_{\frac{1}{4}}(0,1)$ and
$0<r\leq \frac{1}{4}$, let $\eta\in C_0^\infty(P_{2r}(x,t))$ be such that
$0\le\eta\le 1$, $\eta\equiv 1$ on $P_r(x,t)$, $|\nabla\eta|\le 4r^{-1}$,
and $|\partial_t\eta|\le 16r^{-2}$.  Let $(v,P^1):\mathbb R^n\times (0,1)\to\mathbb R^n\times \mathbb R$  solve
\begin{equation}\label{stokes}
\begin{cases}
\partial_t v-\Delta v+\nabla P^1=-\nabla\cdot\Big(\eta^2 (u\otimes u+\nabla d\otimes\nabla d-\frac12{|\nabla d|^2}\mathbb{I}_n)\Big) &\ {\rm{in}}\ \mathbb R^n\times (0,1)\\
\ \ \ \ \ \  \ \ \ \ \ \ \ \ \   \nabla\cdot v= 0 &\ {\rm{in}}\ \mathbb R^n\times (0,1)\\
\ \ \ \ \  \ \ \ \ \ \ \ \ \ \ \ \ \ \ \ v=0 &\ {\rm{on}}\ \mathbb R^n\times \{0\}.
\end{cases}
\end{equation}
Define $w: P_r(x,t)\rightarrow \R^n$ by $w=u-v$. Then $w$ solves
the Stokes equation in $P_r(x,t)$:
\beq\label{unqlc3.17} \left\{
\begin{split}
\partial_t w-\Delta w+\nabla Q^1&=0 \quad\mbox{in } P_r(x,t)\\
\nabla\cdot w&=0\quad\mbox{in } P_r(x,t).
\end{split}\right.
\eeq
By the standard theory of linear Stokes' equations, we have that
$w\in C^\infty(P_r(x,t))$ and,
for any $\theta\in (0, 1)$,
\begin{equation}\label{stokes_est}
\|w\|_{L^p_tL^q_x(P_{\theta r}(x,t))}
\le C\theta\ \|w\|_{L^p_tL^q_x(P_r(x,t))}.
\end{equation}
To estimate $v$, we apply $\mathbb P$ to both sides of the equation (\ref{stokes})$_1$ to obtain
$$\partial_t v-\mathbb A v=-\mathbb P\nabla\cdot\Big(\eta^2 (u\otimes u+\nabla d\otimes\nabla d-\frac12{|\nabla d|^2}\mathbb{I}_n)\Big)
\ {\rm{in}}\ \mathbb R^n\times (0,1); \
v=0 \ {\rm{on}}\ \mathbb R^n\times \{0\}.$$
By the Duhamel formula, we have
\begin{equation}\label{duhamel}
v(t)=-\int_0^te^{-(t-\tau)\mathbb{A}}\mathbb{P}\nabla\cdot\Big(\eta^2(u\otimes u+\nabla d\otimes\nabla d-\frac12{|\nabla d|^2}\mathbb{I}_n)\Big)\,d\tau,
\ 0<t\le 1.
\end{equation}
Now we can apply Fabes-Jones-Riviere \cite{FJR} Theorem 3.1 (see also Kato \cite{kato} page 474, ($2.3'$)) to conclude
that $v\in L^p_tL^q_x(\mathbb R^n\times [0,1])$ and
\begin{eqnarray}\label{v-estimate}
\|v\|_{L^p_tL^q_x(\mathbb R^n\times [0,1])}
&\le& C(\|\eta u\|_{L^p_tL^q_x(\mathbb R^n\times [0,1])}^2+\|\eta \nabla d\|_{L^p_tL^q_x(\mathbb R^n\times [0,1])}^2)\nonumber\\
&\le& C\epsilon_0 (\|u\|_{L^p_tL^q_x(P_{2r}(x,t)))}+\|\nabla d\|_{L^p_tL^q_x(P_{2r}(x,t))}).
\end{eqnarray}
Putting (\ref{stokes_est}) and (\ref{v-estimate}) together, we have that for any $\theta\in (0, 1)$,
\begin{equation}\label{decay_est1}
\|u\|_{L^p_tL^q_x(P_{\theta r}(x,t))}
\le C(\theta+\epsilon_0)\|u\|_{L^p_tL^q_x(P_{2r}(x,t))}
+C\epsilon_0\|\nabla d\|_{L^p_tL^q_x(P_{2r}(x,t))}.
\end{equation}
By Claim 1, we have that
for any $\alpha\in (0,1)$, there exists $\epsilon_0>0$ depending on $\alpha$ such that
\begin{equation}\label{Dd_estimate}
\|\nabla d\|_{L^p_tL^q_x(P_{2r}(x,t))}\le Cr^{\alpha}\|\nabla d\|_{L^p_tL^q_x(P_1(0,1))}.
\end{equation}
Substituting (\ref{Dd_estimate}) into (\ref{decay_est1}) yields
\begin{equation}\label{decay_est2}
\|u\|_{L^p_tL^q_x(P_{\theta r}(x,t))}
\le C(\theta+\epsilon_0)\|u\|_{L^p_tL^q_x(P_{2r}(x,t))}+Cr^\alpha\|\nabla d\|_{L^p_tL^q_x(P_{1}(0,1))}.
\end{equation}
It is standard that by choosing $\theta=\theta_0(\alpha)>0$ and
iterating (\ref{decay_est2}) finitely many times, we conclude that
for any $(x,t)\in P_{\frac{1}{4}}$,
$0<r\leq\frac{1}{4}$ and $0<\alpha<1$,
\begin{equation}\label{decay_est3}
\|u\|_{L^p_tL^q_x(P_{r}(x,t))}
\le C\Big(\|u\|_{L^p_tL^q_x(P_1(0,1))}+\|\nabla d\|_{L^p_tL^q_x(P_{1}(0,1))}\Big)r^\alpha.
\end{equation}
By H\"older's inequality, (\ref{decay_est3}) implies that
$u\in \mathcal M^{2,2-2\alpha}(P_\frac38(0,1))$, and
\begin{equation}\label{u_morrey_estimate}
\|u\|_{\mathcal M^{2,2-2\alpha}(P_\frac38(0,1))}
\le C\Big[\|u\|_{L^p_tL^q_x(P_1(0,1))}+\|\nabla d\|_{L^p_tL^q_x(P_1(0,1))}\Big].
\end{equation}

The higher integrability estimate of $u$ on $P_\frac14(0,1)$ can be done by the parabolic Riesz potential estimate
in parabolic Morrey spaces. Here we will sketch it.
Let $\phi\in C_{0}^{\infty}(P_{\frac{3}{8}}(0,1))$ such that $0\leq\phi\leq 1$, $\phi\equiv1$ on
$P_{\frac{5}{16}}(0,1)$, and
$$|\partial_t\phi|+|\nabla \phi|+|\nabla^2\phi|\leq C.$$
Define $\widetilde{u}:\mathbb R^n\times [0,1]\to\mathbb R^n$ by
\begin{equation}\label{duhamel1}
\widetilde{u}(t)=-\int_0^te^{-(t-\tau)\mathbb{A}}\mathbb{P}\nabla\cdot\Big(\phi^2(u\otimes u+\nabla d\otimes\nabla d-\frac12{|\nabla d|^2}\mathbb{I}_n)\Big)\,d\tau,
\ 0<t\le 1.
\end{equation}
Then, as in the proof of Theorem 3.1 (i) of \cite{FJR}, we have that for any $(x,t)\in\mathbb R^n\times (0,1]$,
\begin{equation}\label{duhamel2}
|\widetilde{u}(x,t)|
\le C\int_0^t\int_{\mathbb R^n}\frac{1}{\delta^{n+1}((x,t), (y,s))}
(|\phi u|^2+|\phi\nabla d|^2)(y,s)\,dyds.
\end{equation}
Recall the parabolic Riesz potential of
order $1$, $I_1(\cdot)$, is defined by
$$I_{1}(f)(z):=\int_{\R^{n+1}}\frac{|f(w)|}{\delta^{n+1}(z,w)}\, dw, \ f\in L^1(\mathbb R^{n+1}).$$
Then we have
\begin{equation}\label{para_morrey}
|\widetilde{u}(x,t)|
\le CI_1({F})(x,t), \ (x,t)\in\mathbb R^n\times (0,1],
\end{equation}
where
$${F}=\phi^2(|u|^2+|\nabla d|^2).$$
By H\"older's inequality, (\ref{Dd_morrey_estimate}), and (\ref{u_morrey_estimate}), we have
that ${F}\in \mathcal{M}^{1,2-2\alpha}(\R^{n+1})$ and
\begin{equation}\label{f_estimate}
\|{F}\|_{\mathcal M^{1,2-2\alpha}(\mathbb R^{n+1})}
\le C\Big(\|\nabla d\|_{L^p_tL^q_x(P_1(0,1))}^2+\|u\|_{L^p_tL^q_x(P_1(0,1))}^2\Big).
\end{equation}
Hence, by \cite{huang-wang} Theorem 3.1 (ii), we conclude that
$\widetilde{u}\in \mathcal M_*^{\frac{2-2\alpha}{1-2\alpha},2-2\alpha}(\R^{n}\times [0,1])$,
and
\begin{eqnarray}\label{morrey_estimate_u}
\|\widetilde{u}\|_{\mathcal M_*^{\frac{2-2\alpha}{1-2\alpha},2-2\alpha}(\R^{n}\times [0,1])}
&\le& C\|{F}\|_{\mathcal{M}^{1,2-2\alpha}(\R^{n+1})}\nonumber\\
&\le& C\left(\|\nabla d\|_{L^p_tL^q_x(P_1(0,1))}^2+\|u\|_{L^p_tL^q_x(P_1(0,1))}^2\right).
\end{eqnarray}
As $\lim\limits_{\alpha\uparrow \frac{1}{2}}\frac{2-2\alpha}{1-2\alpha}=+\infty$,
we have that $\widetilde{u}\in L^{\gamma}(P_{\frac{5}{16}}(0,1))$ for any
$1<\gamma<+\infty$, and
\begin{equation}\label{integral_estimate}
\|\widetilde u\|_{L^\gamma(P_\frac{5}{16})}\le C(\gamma)\Big(\|\nabla d\|_{L^p_tL^q_x(P_1(0,1))}^2+\|u\|_{L^p_tL^q_x(P_1(0,1))}^2\Big).
\end{equation}
Set $\widetilde w=u-\widetilde u$ on $P_\frac{5}{16}(0.1)$. Then it follows from (\ref{lce}) and (\ref{duhamel1})
that
$$\partial_t \widetilde w-\Delta\widetilde w+\nabla \widetilde{Q}=0;
\ \nabla\cdot\widetilde w=0\ \  \ {\rm{in}} \ \ P_{\frac{5}{16}}(0,1).$$
By the standard theory of linear Stokes' equations, we have that
$\widetilde w\in L^\infty(P_\frac14(0,1))$, and
\begin{eqnarray}\label{integral_estimate1}
\|\widetilde w\|_{L^\infty(P_\frac14(0,1))}&\le&
C\|\widetilde w\|_{L^1(P_\frac5{16}(0,1))}
\le C\Big(\|u\|_{L^1(P_\frac5{16}(0,1))}+\|\widetilde u\|_{L^1(P_\frac5{16}(0,1))}\Big)\nonumber\\
&\le& C\Big(\|\nabla d\|_{L^p_tL^q_x(P_1(0,1))}+\|u\|_{L^p_tL^q_x(P_1(0,1))}\Big).
\end{eqnarray}
It is clear that (\ref{integral_estimate_u}) follows from (\ref{integral_estimate}) and (\ref{integral_estimate1}).
This completes the proof of Claim 2.

Finally, it is not hard to see that by the $W^{2,1}_\gamma$-theory for the heat equation and the linear Stokes equation,
and the Sobolev embedding theorem, we have that $(u,\nabla d)\in L^\infty(P_\frac18(0,1))$. Then the Schauder's theory
and the bootstrap argument can imply that $(u,d)\in C^\infty(P_\frac1{16}(0,1))$.
Furthermore,  the estimate (\ref{unqlc3.4}) holds.
This completes the proof. \endpf

\vspace{2mm}
\noindent{\bf Proof of Corollary \ref{reg-th-lcf}}: It is easy to see that when $p> 2$, $q>n$, for any $(x,t)\in\mathbb R^n\times(0,T]$, we can find $R_0>0$ such that
\beq\label{u-lcf-3}
\|u\|_{L^p_t L^q_x(P_{R_0}(x,t))}
+\|\nabla d\|_{L^p_t L^q_x(P_{R_0}(x,t))}\leq \epsilon_0,
\eeq
where $\epsilon_0$ is given in Lemma \ref{unqlc-lemma3.1}. By Theorem \ref{loc-reg-th-lcf}, we conclude that
$(u,d)\in C^{\infty}(P_{\frac{R_0}{16}}(x,t))$. This completes the proof of Theorem \ref{reg-th-lcf}
\endpf

\section{Proof of Theorem \ref{unique3}}
\setcounter{equation}{0}
\setcounter{theorem}{0}

\medskip

In this section, we will prove Theorem \ref{unique3}.
To do this, we need the
following estimate.

\blm\label{unqlc-lemma3.2}{For $T>0$, suppose that $(u,d)$ is a
weak solution to (\ref{lce}) in $\mathbb R^n\times (0,T]$, which
satisfies the assumption of Theorem \ref{unique3}. Then $(u,d)\in
C^{\infty}(\R^n\times(0,T],\mathbb R^n\times S^2)$,  and there exists $t_0>0$ such that
for $0<t\le t_0$, it holds
 \beq\label{unqlc3.32} \sup\limits_{0<\tau\leq
t}\sqrt{\tau}\Big(\|u(\tau)\|_{L^{\infty}(\R^n)}+\|\nabla
d(\tau)\|_{L^{\infty}(\R^n)}\Big)\leq C\Big(\|u\|_{L^p_t L^q_x(\R^n\times [0,t])}
+\|\nabla d\|_{L^p_t L^q_x(\R^n\times [0,t])}\Big).
 \eeq
 In particular, we have
\beq\label{unqlc3.36}
\lim\limits_{t\downarrow
0^+}\sqrt{t}\Big(\|u\|_{L^{\infty}(\R^n)}+\|\nabla
d\|_{L^{\infty}(\R^n)}\Big)=0. \eeq}
\elm

\pf Let $\epsilon_0$ be given by Lemma \ref{unqlc-lemma3.1}.
Since $p>2$ and $q>n$ satisfy (\ref{serrin-condition}),
for any $0<\epsilon\leq\epsilon_0$ we can find $t_0>0$ such that for any $0<\tau\leq \sqrt{t_0}$
\beq\label{u-lcf-3}
\|u\|_{L^p_t L^q_x(\mathbb R^n\times[0,\tau^2])}
+\|\nabla d\|_{L^p_t L^q_x(\mathbb R^n\times[0,\tau^2])}\leq \epsilon.
\eeq
For any $x_0\in \mathbb R^n$, define
 \beq\notag
\begin{split}
&\bar{u}(y,s)=\tau u(x_0+y\tau ,s\tau ^2)\\
&\bar{P}(y,s)=\tau^2P(x_0+y\tau ,s\tau ^2)\\
&\bar{d}(y,s)=d(x_0+y\tau ,s\tau ^2).
\end{split}
 \eeq
  Then $(\bar{u},\bar{P},\bar{d})$ is a weak solution to
(\ref{lce}) on $P_1(0,1)$, and by (\ref{u-lcf-3}),
 \beq\label{u-lcf-4} \|\bar u\|_{L^p_t
L^q_x(P_1(0,1))} +\|\nabla \bar d\|_{L^p_t L^q_x(P_1(0,1))}\leq
\epsilon.
 \eeq
  By Lemma \ref{unqlc-lemma3.1}, we conclude that
 \beq\label{u-lcf-5}
  |\bar u(0,1)|+|\nabla \bar d(0,1)|\leq C
\left(\|\bar u\|_{L^p_t L^q_x(P_1(0,1))} +\|\nabla \bar d\|_{L^p_t
L^q_x(P_1(0,1))}\right).
 \eeq
 By rescaling, this implies
 \beq\label{u-lcf-6}
  \tau\left(| u(x_0,\tau^2)|+|\nabla d(x_0,\tau^2)|\right)\leq C
\left(\|u\|_{L^p_t L^q_x(\mathbb R^n\times[0,\tau^2])}
+\|\nabla d\|_{L^p_t L^q_x(\mathbb R^n\times[0,\tau^2])}\right)\leq C\epsilon.
\eeq
Taking supremum over all $x_0\in\mathbb R^n$ completes the proof.
\endpf

\medskip
\noindent{\bf Proof of Theorem \ref{unique3}}: By (\ref{unqlc3.36}), we have that for any $\epsilon>0$,
there exists $t_0=t_0(\epsilon)>0$ such that
\begin{eqnarray}\label{max_bound}
\mathcal {A}(t_0)&=&\sum_{i=1}^2\Big[\sup_{0\le t\le t_0}\sqrt{t}(\|u_i(t)\|_{L^{\infty}(\R^n)}+\|\nabla
d_i(t)\|_{L^{\infty}(\R^n)})\nonumber\\
&&+(\|u_i\|_{L^p_t L^{q}_x(\R^n\times [0,t_0]))}+\|\nabla
 d_i\|_{L^p_t L^{q}_x(\R^n\times [0,t_0]))})\Big]\le \epsilon.
\end{eqnarray}
It suffices to show $(u_1,d_1)=(u_2,d_2)$ on $\mathbb R^n\times [0,t_0]$. To do so,
let $u=u_1-u_2$ and $d=d_1-d_2$.
Applying $\mathbb{P}$ to both (\ref{lce})$_1$ for $u_1$ and $u_2$ and taking the difference of
resulting equations, we have that
\begin{equation}
\begin{cases}
u_t-\mathbb{A} u=-\mathbb{P}\nabla\cdot\left(u\otimes u_1+u_2\otimes u
+\nabla d\otimes\nabla d_1+\nabla d_2\otimes\nabla d+(|\nabla d_1|+|\nabla d_2|)|\nabla d|\mathbb I_n\right), \\
\nabla\cdot u=0, \\
d_t-\Delta d=[(\nabla d_1+\nabla d_2)\cdot\nabla d\ d_2+|\nabla d_1|^2d]
-[u\cdot\nabla d_1+u_2\cdot\nabla d], \\
(u,d)\Big|_{t=0}=(0,0).
\end{cases}
\end{equation}
By the Duhamel formula, we have that for any $0<t\le t_0$,
 \beq\notag
\begin{split}
u(t)=-\int_0^te^{-(t-\tau)\mathbb{A}}\mathbb{P}\nabla\cdot\Big(u\otimes u_1+u_2\otimes u
+\nabla d\otimes\nabla d_1+\nabla d_2\otimes\nabla d+(|\nabla d_1|+|\nabla d_2|)|\nabla d|\mathbb I_n\Big)\,d\tau,
\end{split}
\eeq
\beq\label{unqlc4.6}
\begin{split}
d(t)=\int_0^te^{-(t-\tau)\Delta}\Big((\nabla d_1+\nabla d_2)\cdot\nabla
d \ d_2+|\nabla d_1|^2d -u\cdot\nabla d_1-u_2\cdot\nabla
d\Big)\,d\tau.
\end{split}
\eeq
For $0<t\le t_0$, set

$$\Phi(t)=\|u\|_{L^p_tL^{q}_x(\R^n\times [0,t]))}+\|\nabla
 d\|_{L^p_tL^{q}_x(\R^n\times [0,t]))}+\sup\limits_{0\leq\tau\leq t}\|d(\cdot,\tau)\|_{L^{\infty}(\R^n)}.$$
By (\ref{unqlc4.6}) and the standard estimate on the heat kernel,  we obtain that
 \beq\label{unqlc4.9}
\begin{split}
\Big\|\nabla
d(t)\Big\|_{L^{q}(\R^n)}\leq &C\Big[\sum_{i=1}^2\int_0^{t}(t-\tau)^{\frac{1}{p}-1}
\|\nabla d_i\|_{L^{q}(\R^n)}\|\nabla
d\|_{L^{q}(\R^n)}\,d\tau\\
&+\|d\|_{L^{\infty}(\R^n)}\int_0^{t}(t-\tau)^{\frac{1}{p}-1}
\|\nabla d_1\|_{L^{q}(\R^n)}^2\,d\tau\\
&+\int_0^{t}(t-\tau)^{\frac{1}{p}-1}
\|\nabla d_1\|_{L^{q}(\R^n)}\|u\|_{L^{q}(\R^n)}\,d\tau\\
&+\int_0^{t}(t-\tau)^{\frac{1}{p}-1}
\|u_2\|_{L^{q}(\R^n)}\|\nabla d\|_{L^{q}(\R^n)}\,d\tau\Big].
\end{split}
\eeq
 By the standard Riesz potential estimate in $L^p$-spaces (see \cite{FJR} Theorem 3.0), we see that
$\nabla d\in L^p_tL^q_x(\mathbb R^n\times [0,t_0])$, and
 \beq\label{unqlc4.10}
\begin{split}
\Big\|\nabla d\Big\|_{L^p_tL^q_x(\R^n\times [0,t_0])}\leq&C\Big[\sum_{i=1}^2\|\nabla
d_i\|_{L^p_tL^{q}_x(\R^n\times [0,t_0]))}\|\nabla
d\|_{L^p_tL^{q}_x(\R^n\times [0,t_0])}\\
&+\|d\|_{L^{\infty}(\R^n\times[0,t_0])}\|\nabla
d_1\|^2_{L^p_tL^{q}_x(\R^n\times [0,t_0])}\\
&+\|\nabla
d_1\|_{L^p_tL^{q}_x(\R^n\times [0,t_0])}\|u\|_{L^p_tL^{q}_x(\R^n\times [0,t_0])}\\
&+\|u_2\|_{L^p_tL^{q}_x(\R^n\times [0,t_0])}\|\nabla
d\|_{L^p_tL^{q}_x(\R^n\times [0,t_0])}\Big]\\
\leq&C\mathcal A(t_0)\Phi(t_0).
\end{split}
\eeq
 Similarly, by using the estimate of Theorem 3.1 (i) of \cite{FJR},  we have that
$u\in L^p_tL^q_x(\mathbb R^n\times [0,t_0])$,  and
 \beq\label{unqlc4.11}
\begin{split}
\|u\|_{L^p_tL^{q}_x(\R^n\times [0,t_0])} \leq& C\mathcal A(t_0)\Phi(t_0).
\end{split}
\eeq
Now we need to  estimate
$\sup\limits_{0\leq\tau\leq t_0}\|d(\cdot,\tau)\|_{L^{\infty}(\R^n)}$. We claim
\begin{equation}\label{max_estimate_d}
\|d\|_{L^\infty(\mathbb R^n\times [0,t_0])}\le C\mathcal A(t_0)\Phi(t_0).
\end{equation}
To show (\ref{max_estimate_d}), let
$\displaystyle H(x,t)$
be the heat kernel of $\R^n$.
By (\ref{unqlc4.6}), we have
 \beq\label{unqlc3.35}
\begin{split}|d(x,t)|=&\Big|\int_0^t\int_{\R^n}H(x-y,t-\tau)\left((\nabla d_1+\nabla d_2)\cdot\nabla d\ d_2+|\nabla
d_1|^2d\right)(y,\tau)\,dyd\tau\\
&-\int_0^t\int_{\R^n}H(x-y,t-\tau)\left(u\cdot\nabla d_1+u_2\cdot\nabla d\right)(y,\tau)\,dyd\tau\Big|\\
\leq& C\Big[\int_0^t\int_{\mathbb R^n} H(x-y,t-\tau)K(y,\tau)\,dyd\tau\\
&\quad +\int_0^t\int_{\R^n} H(x-y,t-\tau)|\nabla d_1|^2(y,\tau)\,dyd\tau
 \cdot\sup\limits_{0\leq\tau\leq
 t}\|d(\cdot,\tau)\|_{L^{\infty}(\R^n)}\Big], \\
\end{split}\eeq
 where
$$K(y,\tau):=\sum\limits_{i=1}^2(|u_i|+|\nabla d_i|)(|u|+|\nabla d|)(y,\tau).$$
 By (\ref{max_bound}), we have that for any $0<t\le t_0$,
 \beq\label{unqlc3.37}
\begin{split}
&\int_0^t\int_{\mathbb R^n} H(x-y,t-\tau)K(y,\tau)\,dyd\tau\\
\leq&\mathcal{A}(t_0)\int_0^{t}(t-\tau)^{-\frac{n}2}\tau^{-\frac12}\int_{\mathbb R^n}(|u|+|\nabla d|)\exp\Big(-\frac{|x-y|^2}{4(t-\tau)}\Big)
\,dyd\tau\\
\leq&\mathcal{A}(t_0)\Big\|(t-\tau)^{-\frac{n}{2q}}\tau^{-\frac12}\Big\|_{L^{\frac{p}{p-1}}([0,t])}
\Big\||u|+|\nabla d|\Big\|_{L^p_tL^q_x(\mathbb R^n\times [0,t])}\\
\leq& C\mathcal {A}(t_0)\Phi(t_0),
\end{split}\eeq
where we have used H\"older inequality and
\begin{eqnarray*}
\Big\|(t-\tau)^{-\frac{n}{2q}}\tau^{-\frac12}\Big\|_{L^{\frac{p}{p-1}}([0,t])}^{\frac{p}{p-1}}
&=&t^{(\frac12-(\frac{n}{2q}+\frac{1}{p}))\frac{p}{p-1}}
\int_0^1 (1-\tau)^{-\frac{np}{2(p-1)q}}\tau^{-\frac{p}{2(p-1)}}\,d\tau\\
&=&\int_0^1 (1-\tau)^{-\frac{p-2}{2(p-1)}}\tau^{-\frac{p}{2(p-1)}}\,d\tau<+\infty,
\end{eqnarray*}
as (i) $\frac{n}{2q}+\frac{1}{p}=\frac12$,
and (ii) $2<p<+\infty$ yields $\frac{p}{2(p-1)}<1$ and $\frac{p-2}{2(p-1)}<1$.

Similarly, we can obtain that for $0\le t\le t_0$,
\begin{equation}\label{unqlc3.38}
\int_0^t\int_{\R^n} H(x-y,t-\tau)|\nabla d_1|^2(y,\tau)\,dyd\tau
\le C\mathcal {A}^2(t_0).
\end{equation}
Putting (\ref{unqlc3.37}) and (\ref{unqlc3.38}) into (\ref{unqlc3.35}) and taking supremum over $(x,t)\in\mathbb R^n\times [0,t_0]$, we have
\begin{equation}\label{max_bound1}
\sup_{0\le t\le t_0}\|d\|_{L^\infty(\mathbb R^n)}
\le C\mathcal{A}(t_0)\Phi(t_0)+C\mathcal {A}^2(t_0)\sup_{0\le t\le t_0}\|d\|_{L^\infty(\mathbb R^n)}.
\end{equation}
Therefore, if we choose $\epsilon\le \sqrt{\frac{1}{2C}}$ so that
$C\mathcal{A}^2(t_0)\le C\epsilon^2\le \frac12$, then we obtain (\ref{max_estimate_d}).

Putting (\ref{unqlc4.10}), (\ref{unqlc4.11}) and (\ref{max_estimate_d}) together, and choosing $\epsilon\le\frac{1}{2C}$,
we obtain
$$\Phi(t_0)\leq C\mathcal{A}(t_0)\Phi(t_0)\le \frac12\Phi(t_0).$$
This implies that $\Phi(t_0)=0$ and hence $(u_1,d_1)\equiv (u_2,d_2)$ on $\mathbb R^n\times [0,t_0]$.
If $t_0<T$, then we can repeat the argument for $t\in [t_0,T]$ and eventually show that $(u_1,d_1)\equiv (u_2,d_2)$ on $\mathbb R^n\times [0,T]$.
This completes the proof. \endpf

\vspace{4mm}
\noindent{\bf Acknowledgements}. The paper is part of my Ph.D. thesis in University of Kentucky. I would like to thank my advisor Professor Changyou Wang for his helpful discussion and constant encouragement.

\end{document}